\documentclass[11pt]{article}
\usepackage{amsfonts}
\usepackage{amsmath}
\usepackage{amssymb}
\usepackage{graphicx}
\usepackage[dvipsnames,usenames]{color}
\usepackage[pdfstartview=FitH ]{hyperref}
\usepackage[hyperpageref]{backref}
\begin{document}
\baselineskip 18pt
\def\o{\over}
\def\e{\varepsilon}
\title{\Large\bf On\ \ the\ \ Inverse\ \ Problem\ \ Relative\ \ to \\
Dynamics\ \ of\ \ the\ \ $w$\ \ Function}
\author{Chaohua\ \  Jia}
\date{}
\maketitle {\small \noindent {\bf Abstract.} In this paper we shall
study the inverse problem relative to dynamics of the w function
which is a special arithmetic function and shall get some results.}

\vskip.5in
\noindent{\bf 1. Introduction}

In 2006, Wushi Goldring [4] proposed some problems and conjectures
on dynamics of the $w$ function and gave some interesting results.
Recently Yong-Gao Chen and Ying Shi [2], [3] made further progress
on these problems. In this paper we shall study the inverse problem
relative to dynamics of the $w$ function.

We begin by introducing some notations. Let ${\cal P}$ be the set of
prime numbers and $P(n)$ denote the largest prime factor of integer
$n>1$. Write
\begin{align*}
C_3&=\{p_1p_2p_3:\ p_i\in {\cal P}\,(i=1,2,3),\ p_i\ne p_j\,(i\ne j)
\},\\
B_3&=\{p_1p_2p_3:\ p_i\in {\cal P}\,(i=1,2,3),\ p_1=p_2\ {\rm or}\
p_1=p_3\\
&\qquad\qquad{\rm or}\ p_2=p_3,\ \hbox{but not}\ p_1=p_2=p_3\},\\
D_3&=\{p^3:\ p\in {\cal P}\}.
\end{align*}
Then
$$
\{p_1p_2p_3:\ p_i\in {\cal P}\,(i=1,2,3)\}=C_3\cup B_3\cup D_3,
$$
where no any two of $C_3,\, B_3$ and $D_3$ intersect. Let
$$
A_3=C_3\cup B_3.
$$
For $n=p_1p_2p_3\in A_3$, define the $w$ function by
$$
w(n)=P(p_1+p_2)P(p_1+p_3)P(p_2+p_3)
$$
and define
$$
w^0(n)=n,\ \ w^i(n)=w(w^{i-1}(n)),\ \ i=1,\,2,\,\cdots,
$$
which is reasonable according to Lemma 1 in section 2.

Wushi Goldring [4] proved that for any $n\in A_3$, there exists $i$
such that $w^i(n)=20$. The smallest of such $i$ is denoted by ${\rm
ind}(n)$. He [4] proposed two conjectures on ${\rm ind}(n)$
(Conjectures 2.9 and 2.10) and gave the first upper bound for ${\rm
ind}(n)$ which has been improved greatly by Yong-Gao Chen and Ying
Shi [2] recently. Wushi Goldring [4] also asked the following
inverse problems:

1. For $n\in A_3$, can we find $m\in A_3$ such that $w(m)=n$?

2. If so, how many such elements are there?

3. What form do they have?

For $n\in A_3$, if there is $m\in A_3$ such that $w(m)=n$, then we
call $m$ a parent of $n$. If this $m\in S\subset A_3$, then we call
it $S$-parent of $n$. Wushi Goldring [4] proved that there are
infinitely many elements of $B_3$ which have at least seven parents.
He proposed the following conjecture (Conjecture 2.16 in [4]).

{\bf Conjecture (Wushi Goldring).} Every element of $A_3$
(respectively $B_3$) has infinitely many $C_3$-parents (respectively
$B_3$-parents).

Yong-Gao Chen and Ying Shi [3] proved that for any given positive
integer $k$, there are infinitely many elements of $B_3$ which have
at least $k$ $B_3$-parents. On the other hand, they [3] proved that
there are infinitely many elements of $B_3$ which have no
$B_3$-parent.

In this paper we shall study parents of elements of $C_3$. It is
obvious that the element of $C_3$ has no $B_3$-parent. We shall
prove that there are infinitely many elements of $C_3$ which have
enough $C_3$-parents.

In the following, $p,\,p_1,\,p_2,\,p_3,\,q,\,\,r,\,r_1,\,r_2$ denote
prime numbers and $c_1,$ $c_2,\,\cdots$ denote positive constants.
The expression $f\ll g$ means $f=O(g)$. We suppose that $x$ is
sufficiently large throughout.

{\bf Theorem 1.} There exists an element $r_1r_2q$ of $C_3$ which
satisfies $x^{1\o 2}\log x< r_i\leq 2x^{1\o 2}\log x\,(i=1,2),\
q\leq 4x$ and has at least $c_1{x\o \log^4x}$ different
$C_3$-parents $p_1p_2p_3$ with $x< p_i\leq 2x\,(i=1,2,3)$.

We shall also prove that there are infinitely many elements of $B_3$
which have enough $C_3$-parents.

{\bf Theorem 2.} There exists an element $qr^2$ of $B_3$ which
satisfies $q\leq 4x,\ x^{1\o 2}\log x< r\leq 2x^{1\o 2}\log x$ and
has at least $c_2{x\o \log^4 x}$ different $C_3$-parents $p_1p_2p_3$
with $x< p_i\leq 2x\,(i=1,2,3)$.

Moreover, we shall prove that there are infinitely many elements of
$B_3$ that have enough $B_3$-parents, which is a quantitative
improvement on the result of Yong-Gao Chen and Ying Shi [3].

{\bf Theorem 3.} There exists an element $qr^2$ of $B_3$ which
satisfies $x< q\leq 2x,\ x^{1\o 2}\log x< r\leq 2x^{1\o 2}\log x$
and has at least $c_3{x^{1\o 2}\o \log^2 x}$ different $B_3$-parents
$pq^2$ with $x< p\leq 2x$.

\vskip.3in
\noindent{\bf 2. Lemmas}

{\bf Lemma 1.} If $n\in A_3$, then $w(n)\in A_3$.

This is Lemma 2.1 in [4].

{\bf Lemma 2.} Let $n_j\,(1\leq j\leq Z)$ be distinct positive
integers not exceeding $N$ and $Z(N;\,r,\,a)$ denote the number of
those $n_j$ which are congruent to $a\,({\rm mod}\,r)$. If $X\geq
2$, then we have
$$
\sum_{r\leq X}r\sum_{a=1}^r\Bigl(Z(N;\,r,\,a)-{Z\o r}\Bigr)^2\ll
(N+X^2)Z.
$$

This is Theorem 1 in [1], which is obtained by the large sieve
method.

{\bf Lemma 3.} We have
$$
\sum_{x^{1\o 2}\log x< r\leq 2x^{1\o 2}\log x}\sum_{x< p_1\leq 2x}
\Bigl(\sum_{\substack{x< p\leq 2x\\ p\equiv -p_1\,({\rm mod}\,r)}}1
-{1\o r}\sum_{x< p\leq 2x} 1 \Bigr)^2\ll {x^2\o \log x}.
$$

{\bf Proof.} We see
\begin{align*}
&\ \,\sum_{x^{1\o 2}\log x< r\leq 2x^{1\o 2}\log x}\sum_{x< p_1\leq
2x}\Bigl(\sum_{\substack{x< p\leq 2x\\ p\equiv -p_1\,({\rm mod}\,r)}}
1 -{1\o r}\sum_{x< p\leq 2x} 1 \Bigr)^2\\
&\leq\sum_{x^{1\o 2}\log x< r\leq 2x^{1\o 2}\log x}\sum_{x< n_1\leq
x+([{x\o r}]+1)r}\Bigl(\sum_{\substack{x< p\leq 2x\\ p\equiv
-n_1\,({\rm mod}\,r)}}1 -{1\o r}\sum_{x< p\leq 2x} 1 \Bigr)^2
\end{align*}
\begin{align*}
&\ll\sum_{x^{1\o 2}\log x< r\leq 2x^{1\o 2}\log x}{x\o
r}\,\sum_{a=1}^r \Bigl(\sum_{\substack{x< p\leq 2x\\ p\equiv
a\,({\rm mod}\,r)}}1 -{1\o r}\sum_{x< p\leq 2x} 1 \Bigr)^2\\
&\ll{1\o \log^2x}\sum_{r\leq 2x^{1\o 2}\log x}r\,\sum_{a=1}^r
\Bigl(\sum_{\substack{x< p\leq 2x\\ p\equiv a\,({\rm mod}\,r)}}1
-{1\o r}\sum_{x< p\leq 2x} 1 \Bigr)^2.
\end{align*}
Then Lemma 2 and the prime number theorem yield
\begin{align*}
&\ \,{1\o \log^2x}\sum_{r\leq 2x^{1\o 2}\log x}r\,\sum_{a=1}^r
\Bigl(\sum_{\substack{x< p\leq 2x\\ p\equiv a\,({\rm mod}\,r)}}1
-{1\o r}\sum_{x< p\leq 2x} 1 \Bigr)^2\\
&\ll{1\o \log^2x}\cdot x\log^2x\cdot{x\o \log x}\\
&={x^2\o \log x}.
\end{align*}
Hence, Lemma 3 holds true.

\vskip.3in
\noindent{\bf 3. The proof of Theorem 1}

We note that if $p> n^{1\o 2}$, then $p|\,n\Longleftrightarrow P(n)=
p$. Then we have
\begin{align*}
&\ \,\sum_{x< p_1\leq 2x}\sum_{\substack{x< p_2\leq 2x\\ x^{1\o
2}\log x< P(p_1+p_2)\leq 2x^{1\o 2}\log x}}\sum_{\substack{x< p_3
\leq 2x\\ x^{1\o 2}\log x< P(p_1+p_3)\leq 2x^{1\o 2}\log x}} 1\\
&=\sum_{x< p_1\leq 2x}\Bigl(\sum_{\substack{x< p\leq 2x\\ x^{1\o
2}\log x< P(p+p_1)\leq 2x^{1\o 2}\log x}} 1\Bigr)^2\\
&=\sum_{x<p_1\leq 2x}\Bigl(\sum_{x^{1\o 2}\log x< r\leq 2x^{1\o 2}
\log x}\sum_{\substack{x< p\leq 2x\\ P(p+p_1)=r}} 1\Bigr)^2\\
&=\sum_{x< p_1\leq 2x}\Bigl(\sum_{x^{1\o 2}\log x< r\leq 2x^{1\o 2}
\log x}\sum_{\substack{x< p\leq 2x\\ p\equiv -p_1\,({\rm mod}\,r)}}
1\Bigr)^2\\
&=\sum_{x< p_1\leq 2x}\Bigl(\sum_{x^{1\o 2}\log x< r\leq 2x^{1\o 2}
\log x}{1\o r}\sum_{x< p\leq 2x}1\\
&\ +\sum_{x^{1\o 2}\log x< r\leq 2x^{1\o 2}\log
x}\Bigl(\sum_{\substack{x< p\leq 2x\\ p\equiv -p_1\,({\rm mod}\,r)}}
1-{1\o r}\sum_{x< p\leq 2x}1\Bigr)\Bigr)^2 \tag 1
\end{align*}

\begin{align*}
&=\sum_{x< p_1\leq 2x}\Bigl(\sum_{x^{1\o 2}\log x< r\leq 2x^{1\o
2}\log x}{1\o r}\sum_{x< p\leq 2x}1\Bigr)^2\\
&\ +2\sum_{x< p_1\leq 2x}\Bigl(\sum_{x^{1\o 2}\log x< r\leq 2x^{1\o
2}\log x}{1\o r}\sum_{x< p\leq 2x}1\Bigr)\\
&\ \cdot\sum_{x^{1\o 2}\log x< r\leq 2x^{1\o 2}\log
x}\Bigl(\sum_{\substack{x< p\leq 2x\\ p\equiv -p_1\,({\rm mod}\,r)}}
1-{1\o r}\sum_{x< p\leq 2x}1\Bigr)\\
&+\sum_{x< p_1\leq 2x}\Bigl(\sum_{x^{1\o 2}\log x< r\leq 2x^{1\o 2}
\log x}\Bigl(\sum_{\substack{x< p\leq 2x\\ p\equiv -p_1\,({\rm
mod}\,r)}} 1-{1\o r}\sum_{x< p\leq 2x}1\Bigr)\Bigr)^2.
\end{align*}

The prime number theorem yields
\begin{align*}
&\ \,\sum_{x< p_1\leq 2x}\Bigl(\sum_{x^{1\o 2}\log x<
r\leq 2x^{1\o 2}\log x}{1\o r}\sum_{x< p\leq 2x}1\Bigr)^2\\
&\gg {x^3\o \log^3x}\Bigl(\sum_{x^{1\o 2}\log x< r\leq 2x^{1\o 2}
\log x}{1\o r}\Bigr)^2\gg {x^3\o \log^5x}.
\end{align*}
By the Cauchy inequality and Lemma 3, we have
\begin{align*}
&\ \,\sum_{x< p_1\leq 2x}\Bigl(\sum_{x^{1\o 2}\log x< r\leq 2x^{1\o
2} \log x}\Bigl(\sum_{\substack{x< p\leq 2x\\ p\equiv -p_1\,({\rm
mod}\,r)}} 1-{1\o r}\sum_{x< p\leq 2x}1\Bigr)\Bigr)^2\\
&\leq\sum_{x< p_1\leq 2x}\Bigl(\sum_{x^{1\o 2}\log x< r\leq 2x^{1\o
2} \log x}1\cdot\sum_{x^{1\o 2}\log x< r\leq 2x^{1\o 2}\log
x}\Bigl(\sum_{\substack{x< p\leq 2x\\ p\equiv -p_1\,({\rm mod}\,r)}} 1\\
&\qquad\qquad\qquad\qquad-{1\o r}\sum_{x< p\leq 2x}1\Bigr)^2\Bigr)\\
&\ll x^{1\o 2}\sum_{x^{1\o 2}\log x< r\leq 2x^{1\o 2}\log x}\sum_{x<
p_1\leq 2x}\Bigl(\sum_{\substack{x< p\leq 2x\\ p\equiv -p_1\,({\rm
mod}\,r)}}1-{1\o r}\sum_{x< p\leq 2x}1\Bigr)^2\\
&\ll x^{5\o 2}.
\end{align*}
We also have
$$
\sum_{x< p_1\leq 2x}\Bigl(\sum_{x^{1\o 2}\log x<
r\leq 2x^{1\o 2}\log x}{1\o r}\sum_{x< p\leq 2x}1\Bigr)
$$
\begin{align*}
&\ \cdot\sum_{x^{1\o 2}\log x< r\leq 2x^{1\o 2}\log
x}\Bigl(\sum_{\substack{x< p\leq 2x\\ p\equiv -p_1\,({\rm mod}\,r)}}
1-{1\o r}\sum_{x< p\leq 2x}1\Bigr)\\
&\ll\Bigl(\sum_{x< p_1\leq 2x}\Bigl(\sum_{x^{1\o 2}\log x< r\leq
2x^{1\o 2} \log x}{1\o r}\sum_{x< p\leq 2x}1\Bigr)^2\Bigr)^{1\o 2}\\
&\cdot\Bigl(\sum_{x< p_1\leq 2x}\Bigl(\sum_{x^{1\o 2}\log x< r\leq
2x^{1 \o 2}\log x}\Bigl(\sum_{\substack{x< p\leq 2x\\ p\equiv
-p_1\,({\rm mod}\,r)}}1-{1\o r}\sum_{x< p\leq 2x}1\Bigr)\Bigr)^2\Bigr)
^{1\o 2}\\
&\ll \Bigl(x^3\Bigl(\sum_{x^{1\o 2}\log x< r\leq 2x^{1\o 2}\log
x}{1\o r}
\Bigr)^2\Bigr)^{1\o 2}(x^{5\o 2})^{1\o 2}\\
&\ll x^{3\o 2}\cdot x^{5\o 4}=x^{11\o 4}.
\end{align*}

Combining the above estimates, we get
\begin{align*}
&\ |\{p_1p_2p_3:\ x< p_i\leq 2x\,(i=1,\,2,\,3),\ x^{1\o
2}\log x< P(p_1+p_2)
\leq 2x^{1\o 2}\log x,\\
&\qquad\qquad\qquad x^{1\o 2}\log x< P(p_1+p_3)\leq 2x^{1\o 2}\log x\}|\\
&\gg {x^3\o \log^5x}.
\end{align*}

Now we shall confine $w(p_1p_2p_3)$ to $C_3$. We have
\begin{align*}
&\ \,\sum_{x< p_1\leq 2x}\sum_{\substack{x< p_2\leq 2x\\ x^{1\o
2}\log x< P(p_1+p_2)\leq 2x^{1\o 2}\log x}}\sum_{\substack{x< p_3\leq 2x\\
x^{1\o 2}\log x< P(p_1+p_3)\leq 2x^{1\o 2}\log x\\ P(p_1+p_3)=P(p_2+p_3)}} 1\\
&\leq\sum_{x^{1\o 2}\log x< r\leq 2x^{1\o 2}\log x}\sum_{x< p_3\leq
2x}\sum_{\substack{x< p_1\leq 2x\\ P(p_1+p_3)=r}}\sum_{\substack{x<
p_2\leq 2x\\ P(p_2+p_3)=r}} 1\\
&=\sum_{x^{1\o 2}\log x< r\leq 2x^{1\o 2}\log x}\sum_{x< p_3\leq 2x}
\sum_{\substack{x< p_1\leq 2x\\
p_1\equiv -p_3\,({\rm mod}\,r)}}
\sum_{\substack{x< p_2\leq 2x\\ p_2\equiv -p_3\,({\rm mod}\,r)}} 1\\
&\leq\sum_{x^{1\o 2}\log x< r\leq 2x^{1\o 2}\log x}\sum_{x< n_3\leq
2x}\sum_{\substack{x< n_1\leq 2x\\ n_1\equiv -n_3\,({\rm mod}\,r)}}
\sum_{\substack{x< n_2\leq 2x\\ n_2\equiv -n_3\,({\rm mod}\,r)}}1\\
&\ll\sum_{x< n_3\leq 2x}\sum_{x^{1\o 2}\log x< r\leq 2x^{1\o 2}\log
x}{x^2\o r^2}\\
&\ll x^3\sum_{x^{1\o 2}< n}{1\o n^2}\ll x^{5\o 2}.
\end{align*}
Similarly,
$$
\sum_{x< p_1\leq 2x}\sum_{\substack{x< p_2\leq 2x\\ x^{1\o 2}\log x<
P(p_1+p_2) \leq 2x^{1\o 2}\log x}}\sum_{\substack{x< p_3\leq 2x\\
x^{1\o 2}\log x< P(p_1+p_3)\leq 2x^{1\o 2}\log x\\
w(p_1p_2p_3)\not\in C_3}}1 =O(x^{5\o 2}).
$$
Therefore
\begin{align*}
&\ |\{p_1p_2p_3\in C_3:\ x< p_i\leq 2x\,(i=1,\,2,\,3),\ x^{1\o
2}\log x< P(p_1+p_2)\leq 2x^{1\o 2}\log x,\\
&\qquad\qquad x^{1\o 2}\log x< P(p_1+p_3)\leq 2x^{1\o 2}\log x,\
w(p_1p_2p_3)\in C_3\}| \tag 2\\
&\gg {x^3\o \log^5x}.
\end{align*}

On the other hand, the number of triples $(r_1,\,r_2,\,q)$ is
$O({x^2\o \log x})$, where $x^{1\o 2}\log x< r_i\leq 2x^{1\o 2}\log
x\,(i=1,\,2),\ q\leq 4x$. Therefore in the set in (2), there are at
least $c_4{x\o \log^4 x}$ different triples $(p_1,\,p_2,\,p_3)$
satisfying $P(p_1+p_2)=r_1,\ P(p_1+p_3)=r_2,\ P(p_2+ p_3)=q$ for
some $(r_1,\,r_2,\,q)$. In other words, there are at least ${1\o
3!}c_4{x\o \log^4 x}$ different numbers $n=p_1p_2p_3\in C_3$ such
that $w(n)=$ $r_1r_2q$.

So far the proof of Theorem 1 is finished.

\vskip.3in
\noindent{\bf 4. The proof of Theorem 2}

We have
\begin{align*}
&\ \,\sum_{x< p_1\leq 2x}\sum_{x< p_2\leq 2x}\sum_{\substack{x< p_3
\leq 2x\\ x^{1\o 2}\log x< P(p_1+p_2)=P(p_1+p_3)\leq 2x^{1\o 2} \log
x}}1\\
&=\sum_{x^{1\o 2}\log x< r\leq 2x^{1\o 2}\log x}\sum_{x< p_1\leq 2x}
\sum_{\substack{x< p_2\leq 2x\\ P(p_2+p_1)=r}}\sum_{\substack{x< p_3
\leq 2x\\ P(p_3+p_1)=r}} 1\\
&=\sum_{x^{1\o 2}\log x< r\leq 2x^{1\o 2}\log x}\sum_{x< p_1\leq 2x}
\Bigl(\sum_{\substack{x< p\leq 2x\\ P(p+p_1)=r}} 1\Bigr)^2\\
&=\sum_{x^{1\o 2}\log x< r\leq 2x^{1\o 2}\log x}\sum_{x< p_1\leq 2x}
\Bigl(\sum_{\substack{x< p\leq 2x\\ p\equiv -p_1\,({\rm mod}\,r)}}
1\Bigr)^2.
\end{align*}

By the inequality
$$
a^2+b^2\geq {1\o 2}\,(a-b)^2,
$$
we can get
\begin{align*}
&\ \,\sum_{x^{1\o 2}\log x< r\leq 2x^{1\o 2}\log x}\sum_{x< p_1\leq
2x} \Bigl(\sum_{\substack{x< p\leq 2x\\ p\equiv -p_1\,({\rm
mod}\,r)}} 1\Bigr)^2\\
&\geq \sum_{x^{1\o 2}\log x< r\leq 2x^{1\o 2}\log x}\sum_{x< p_1\leq
2x}{1\o 2}\Bigl({1\o r}\sum_{x< p\leq 2x} 1\Bigr)^2\\
&-\sum_{x^{1\o 2}\log x< r\leq 2x^{1\o 2}\log x}\sum_{x< p_1\leq 2x}
\Bigl(\sum_{\substack{x< p\leq 2x\\ p\equiv -p_1\,({\rm mod}\,r)}} 1
-{1\o r}\sum_{x< p\leq 2x} 1\Bigr)^2.
\end{align*}

The prime number theorem yields
\begin{align*}
&\ \,\sum_{x^{1\o 2}\log x< r\leq 2x^{1\o 2}\log x}\sum_{x< p_1\leq
2x}{1\o 2}\Bigl({1\o r}\sum_{x< p\leq 2x} 1\Bigr)^2\\
&\gg {x^3\o \log^3x}\sum_{x^{1\o 2}\log x< r\leq 2x^{1\o 2}\log x}{1\o r^2}\\
&\gg {x^{5\o 2}\o \log^5x}
\end{align*}
and Lemma 3 yields
$$
\sum_{x^{1\o 2}\log x< r\leq 2x^{1\o 2}\log x}\sum_{x< p_1\leq
2x}\Bigl (\sum_{\substack{x< p\leq 2x\\ p\equiv -p_1\,({\rm
mod}\,r)}} 1-{1\o r} \sum_{x< p\leq 2x} 1\Bigr)^2=O(x^2).
$$

Therefore
\begin{align*}
&\ |\{p_1p_2p_3\in C_3:\ x< p_i\leq
2x\,(i=1,\,2,\,3),\ x^{1\o 2}\log x<
P(p_1+p_2)=\\
&\qquad\qquad\qquad\qquad P(p_1+p_3)\leq 2x^{1\o 2}\log x\}| \tag 3\\
&\gg {x^{5\o 2}\o \log^5x},
\end{align*}
since the contribution from terms with $p_1p_2p_3\not\in C_3$ is
$O(x^2)$.

On the other hand, the number of $(q,\,r)$ with $q\leq 4x,\ x^{1\o
2}\log x< r\leq 2x^{1\o 2}\log x$ is $O({x^{3\o 2}\o \log x})$.
Therefore in the set in (3), there are at least $c_5\,{x\o \log^4
x}$ different triples $(p_1,\,p_2,\,p_3)$ satisfying
$P(p_1+p_2)=P(p_1+p_3)=r,\ P(p_2+ p_3)=q$ for some $(q,\,r)$. In
other words, there are at least ${1\o 3!}\,c_5\,{x\o \log^4 x}$
different numbers $n=p_1p_2p_3\in C_3$ such that $w(n)=qr^2$. By
Lemma 1, we know $q\ne r$.

So far the proof of Theorem 2 is finished.

\vskip.3in
\noindent{\bf 5. The proof of Theorem 3}

By Lemma 3, we have
\begin{align*}
&\ \,\sum_{x^{1\o 2}\log x< r\leq 2x^{1\o 2}\log x}\sum_{x< q\leq
2x} \Bigl(\sum_{\substack{x< p\leq 2x\\ p\equiv -q\,({\rm
mod}\,r)}}1-{1\o r}\sum_{x< p\leq 2x} 1\Bigr)\\
&\ll \Bigl(\sum_{x^{1\o 2}\log x< r\leq 2x^{1\o 2}\log x}\sum_{x<
q\leq 2x} 1\Bigr)^{1\o 2}\Bigl(\sum_{x^{1\o 2}\log x< r\leq 2x^{1\o
2}\log x}\sum_{x< q\leq 2x}\\
&\qquad\qquad\cdot\Bigl(\sum_{\substack{x< p\leq 2x\\ p\equiv
-q\,({\rm mod}\,r)}}1-{1\o r}\sum_{x< p\leq 2x} 1\Bigr)^2\Bigr)^{1
\o 2}\\
&\ll (x^{3\o 2})^{1\o 2}(x^2)^{1\o 2}=x^{7\o 4}.
\end{align*}
Hence,
\begin{align*}
&\ \,\sum_{x^{1\o 2}\log x< r\leq 2x^{1\o 2}\log x}\sum_{x< q\leq
2x}\sum_{\substack{x< p\leq 2x\\ p\ne q\\ P(p+q)=r}}1\\
&=\sum_{x^{1\o 2}\log x< r\leq 2x^{1\o 2}\log x}\sum_{x< q\leq 2x}
\sum_{\substack{x< p\leq 2x\\ P(p+q)=r}} 1+O(x)\\
&=\sum_{x^{1\o 2}\log x< r\leq 2x^{1\o 2}\log x}\sum_{x< q\leq 2x}
\sum_{\substack{x< p\leq 2x\\ p\equiv -q\,({\rm mod}\,r)}} 1+O(x)
\tag 4\\
&=\sum_{x^{1\o 2}\log x< r\leq 2x^{1\o 2}\log x}\sum_{x< q\leq 2x}
{1\o r}\sum_{x< p\leq 2x}1+O(x^{7\o 4})\\
&\geq {1\o 2}\cdot{x^2\o \log^2x}\sum_{x^{1\o 2}\log x< r\leq
2x^{1\o 2}\log x}{1\o r}+O(x^{7\o 4})\\
&\geq{1\o 3}\cdot{x^2\o \log^3x}.
\end{align*}

Therefore there must be one pair $(q,\,r)$ with $x< q\leq 2x,\ x^{1
\o 2}\log x< r\leq 2x^{1\o 2}\log x$ such that there are at least
${1\o 10}\cdot{x^{1\o 2}\o \log^2 x}$ different $p(\ne q)$
satisfying $x<p\leq 2x,\ P(p+q)=r$. Otherwise we should have
\begin{align*}
&\ \,{1\o 3}\cdot{x^2\o \log^3x}\leq\sum_{x^{1\o 2}\log x< r\leq
2x^{1\o 2} \log x}\sum_{x< q\leq 2x}\sum_{\substack{x< p\leq 2x\\
p\ne q\\ P(p+q)=r}}1\\
&\leq {1\o 10}\cdot{x^{1\o 2}\o \log^2 x}\sum_{x^{1\o 2}\log x<
r\leq 2x^{1\o 2}\log x}\sum_{x< q\leq 2x}1\\
&\leq {1\o 4}\cdot{x^2\o \log^3 x},
\end{align*}
which is a contradiction.

For this pair $(q,\,r)$, there are at least ${1\o 10}\cdot{x^{1\o
2}\o \log^2 x}$ different $p(\ne q)$ such that $w(pq^2)=qr^2$. By
Lemma 1, we know $q\ne r$. Hence, Theorem 3 holds true.

\vskip.3in \noindent{\bf Acknowledgements}

In June 2007, Professor Yong-Gao Chen visited Morningside
Mathematical Center of Academia Sinica in Beijing and gave two talks
to introduce his joint work with Ying Shi on dynamics of the $w$
function. I would like to thank Professor Yong-Gao Chen for his
wonderful talks which attract my interest to this topic. I also
thank all my colleagues and friends in the ``ergodic prime number
theorem'' seminar in Morningside Center for helpful discussion.

Results in this paper were reported at the ``combinatorial and
analytic number theory'' seminar in Nanjing Normal University in
September, 2007 and were reported at the fifth Japan-China seminar
on number theory in Osaka in August, 2008.

\pagebreak

\bigskip

\

Institute of  Mathematics, Academia Sinica, Beijing 100190, P. R.
China

E-mail: jiach@math.ac.cn
\end{document}